\documentstyle[12pt]{report}

\newcommand{\vs}{\vspace{.15in}}

\newcommand{\noin}{\noindent}
\newcommand{\n}{\indent}
\newcommand{\prop}{{\bf Proposition} }
\newcommand{\lemma}{{\bf Lemma}  }
\newcommand{\thm}{{\bf Theorem}  }
\newcommand{\proof}{{\sl Proof.}  }
\newcommand{\rem}{{\bf Remark} }

\newcommand{\qed}{{$\Box $} \vs}
\def\de{{\partial}}
\def\n{{\rm norm}}
\def\pr{{\rm pr}}
\def\C{{\bf C}}

\def\tS{{\tilde S}}
\def\tX{{\tilde X}}
\def\G{\Gamma}

\def\cC{{\cal C}}

\def\cM{{\cal M}}

\def\cP{{\cal P}}
\def\cQ{{\cal Q}}

\def\cI{{\cal I}}
\def\Z{{\bf Z}}

\def\Im{{{\rm Im} \, }}

\begin{document}
\begin{center} {\bf FREE $\C_+$-ACTIONS ON $\C^3$ ARE TRANSLATIONS}
\\[3ex] {\bf SHULIM KALIMAN} \footnote{ The author was partially
supported by the NSA grant MDA904-00-1-0016.}\\[3ex] Department of
Mathematics, University of Miami \\ Coral Gables, Florida 33124,
USA \\ (e-mail: kaliman@math.miami.edu )\\[5ex]
\end{center}
ABSTRACT. {\small Let $X'$ be a smooth contractible
three-dimensional affine algebraic variety with a free algebraic
$\C_+$-action on it such that $S =X'//\C_+$ is smooth. We prove
that $X'$ is isomorphic to $S \times \C$ and the action is induced
by a translation on the second factor. As a consequence we show
that any free algebraic $\C_+$-action on $\C^3$ is a translation
in a suitable coordinate system.} \\[3ex]

{\bf 1. Introduction.} In 1968 Rentschler [Re] proved that every
algebraic action of the additive group $\C_+$ of complex numbers
on $\C^2$ is triangular in a suitable polynomial coordinate
system. This implies that a free $\C_+$-action on $\C^2$ (i.e. an
action without fixed points) can be viewed as a translation. In
1984 Bass [Ba] found a $\C_+$-action on $\C^3$ which is not
triangular in any coordinate system. But the question about free
$\C_+$-actions on $\C^3$ remained open (e.g., see [Sn] , [Kr],
[DaFr], [DeFi]). While working on this problem we consider a more
general situation when there is a nontrivial algebraic
$\C_+$-action on a complex three-dimensional affine algebraic
variety $X'$ such that its ring of regular functions is factorial
and $H_3(X')=0$. By a theorem of Zariski [Za] the algebraic
quotient $X'//\C_+$ is isomorphic to an affine surface $S$. Let
$\pi : X' \to S$ be the natural projection. Then there is a curve
$\G \subset S$ such that for $E= \pi^{-1}( \G )$ the variety $X'
\setminus E$ is isomorphic to $( S \setminus \G ) \times \C$ over
$S \setminus \G$. The study of morphism $\pi |_E : E \to \G$ is
central for this paper. As an easy consequence of the Stein
factorization one can show that $\pi |_E = \theta \circ \vartheta$
where $\vartheta : E \to Z$ is a morphism into a curve $Z$ with
general fibers isomorphic to $\C$, and $\theta : Z \to \G$ is a
quasi-finite morphism. A more delicate fact (Proposition 3.2) is
that $\theta$ is, actually, finite and, furthermore, if in
addition $X'$ is smooth and $H_2(X')=0$, then each irreducible
component $Z^1$ of $Z$, such that $\theta |_{Z^1}$ is not
injective, is rational and $\theta (Z^1)$ is a polynomial curve. A
geometrical observation (Proposition 4.2) shows that in the case
of a free $\C_+$-action and a smooth $S$ such a component $Z^1$
cannot exist. If the restriction of $\theta$ to any irreducible
component of $Z$ is injective then $\G$ can be chosen empty
(Theorem 5.4), i.e. $X'$ is isomorphic to $S \times \C$ over $S$.
In combination with  Miyanishi's theorem [Miy80] this yields the
long-expected result. \vs

\thm 1. {\em Every free algebraic $\C_+$-action on $\C^3$ is a
translation in a suitable polynomial coordinate system.} \vs

{\em Acknowledgments.} It is our pleasure to thank P. Bonnet, D.
Daigle, G. Freudenburg, K.-H. Fieseler, and M. Zaidenberg for
useful consultations.\vs

\noin {\bf 2. 
Decomposition $\pi |_E = \theta \circ \vartheta$. \vs }

\lemma 2.1. {\em Let $X'$ be a factorial affine algebraic variety
of dimension 3 with a nontrivial $\C_+$-action on it and $\pi : X'
\to S$ be the quotient morphism. Then $S$ is a factorial affine
algebraic surface, $\pi^{-1}(s_0)$ is either a curve or empty for
any $s_0 \in S$, and $E=\pi^{-1}(\G )$ is a nonempty irreducible
surface for every closed irreducible curve $\G \subset S$.
Furthermore, if $\G =g^* (0)$ for a regular function $g$ on $S$
then $E =(g\circ \pi )^* (0)$.}

\proof By [Za] $S$ is an affine algebraic surface. Suppose that
$E$ is the union of a surface $E'$ and a curve $C$ where $E'$ and
$C$ may be empty. Assume $\pi (E')$ is not dense in $\G$. Consider
a rational function $h$ on $S$ with poles on $\G$ only. Let $e$ be
the product of $h$ and a regular function that does not vanish at
general points of $\G$ but vanishes on $\pi (E')$ with
sufficiently large multiplicity such that $e \circ \pi$ vanishes
on $E' \setminus C$. As $X'$ is factorial, by deleting
singularities in codimension 2 we see that $e \circ \pi$ is a
regular function on $X'$. As it is invariant under the action, $e$
must be regular. Contradiction. Thus $\pi (E')$ is dense in $\G$.
Since $X'$ is factorial $E'$ is the zero fiber of a regular
function $f$ on $X'$. As $f$ does not vanish on a general fiber
$\pi^{-1}(s) \simeq \C$ it is constant on each general fiber.
Hence $f$ is invariant under the action, i.e. $f=g \circ \pi$
where $g$ a regular function on $S$ which implies that $g$ is a
defining function for $\G$ (i.e. $S$ is factorial) and,
furthermore, $E=E' \ne \emptyset$. If for a regular function $f_1$
on $X'$ its zero fiber is an irreducible component of $E$ then
$f_1$ is again invariant under the action, i.e. $f_1=g_1 \circ
\pi$. As $g_1$ vanishes on $\G$ only it must be proportional to
$g$. Thus $E$ is irreducible and we get the last statement as
well. If we assume that $\pi^{-1}(s_0)$ is a surface for $s_0 \in
S$ the similar argument yields a regular function on $S$ that
vanishes at a finite set only which is impossible.
\qed

If we have a nontrivial $\C_+$-action on a three-dimensional
affine algebraic variety $X'$ with a quotient morphism $\pi : X'
\to S=X'//\C_+$ then it is well-known that there is a closed curve
$\G \subset S$ such that for $ \breve{S}=S\setminus \G$ and $E =
\pi^{-1}(\G )$ the variety $X' \setminus E = \pi^{-1}( \breve{S})$
is naturally isomorphic to $\breve{S } \times \C$ where the
projection to the first factor corresponds to $\pi$. Let $\sigma_0
: X' \setminus E \to \breve{S} \times \C$ be an isomorphism over $
\breve{S}$, i.e. $\sigma_0 = ( \pi |_{X' \setminus E}, h_0)$ where
$h_0$ is a regular function on $X' \setminus E$. There is some
freedom in the choice of $h_0$ and, therefore, $\sigma_0$. \vs

\lemma 2.2. {\em Let $X',E, \G ,S$, and $\sigma_0 $ be as before,
$\G =g^{-1}(0)$ for a regular function $g$ on $S$, and $X = S
\times \C$. Then $\sigma_0$ can be chosen so that there exists and
affine modification $\sigma : X' \to X$ such that $\sigma |_{X'
\setminus E} = \sigma_0$.}

\proof Note that $\sigma_0$ extends to a rational map $\sigma
=(\pi , h) : X' \to X$ where $h$ is a rational function on $X'$
whose restriction to $X'\setminus E$ coincides with $h_0$. One can
replace $h_0$ by $h_0 g^n$ (where we treat $g$ as a function on
$X'$) as $\sigma_0$ remains an isomorphism after this replacement.
Choosing $n>0$ sufficiently large we can make $h$ regular on $X'$
whence $\sigma$ becomes a regular birational morphism. \qed

\lemma 2.3. {\em Let $X'$ be a factorial affine algebraic variety
of dimension 3 with a nontrivial $\C_+$-action on it, $\pi : X'
\to S$ be the quotient morphism, $E$ and $\G$ be as before Lemma
2.2.

{\rm (1)} There exist a curve $Z$, a morphism $\vartheta : E \to
Z$ with connected general fibers, and a quasi-finite morphism
$\theta : Z \to \G$ for which $\pi |_{E} = \theta  \circ
\vartheta$. Furthermore, $\theta^{-1} (\G^1)$ is irreducible for
every irreducible component $\G^1$ of $\G$.

{\rm (2)} There exists a finite set $F \subset \G$ such that for $
\G^* = \G \setminus F$, $Z^*=\theta^{-1}(\G^* )$, and $E^*
=\pi^{-1} (\G^*)$ there is an isomorphism $E^* \simeq Z^* \times
\C$ such that under this isomorphism $\vartheta |_{E^*}$ coincides
with the projection to the first factor.}

\proof Extend $\pi |_{E} : E \to \G$ to a proper morphism $ \bar{E
} \to \bar{\G }$ of complete varieties. By the Stein factorization
theorem this extended morphism is a composition of morphisms $
\bar{E} \to \bar{Z}$ and $ \bar{Z } \to \bar{\G }$ where the first
morphism has connected fibers, $ \bar{Z}$ is a curve, and the last
morphism is finite. The restriction of these morphisms is what we
need in (1). Since $E^1 =\pi^{-1}(\G^1)$ is irreducible by Lemma
2.1, we see that $\theta^{-1}(\G^1)=\vartheta (E^1)$ is also
irreducible. As $E$ is the exceptional divisor of affine
modification $\sigma$ from Lemma 2.2, (2) follows from [Ka02, Cor.
2.2] (in the case of a free $\C_+$-action it follows also from the
fact that the general fibers of $\vartheta$ are orbits of the
action, and, therefore, isomorphic to $\C$). \qed

\rem 2.4. (1) Let the center $C$ of the affine modification
$\sigma$ from Lemma 2.2 (i.e. $C$ is the closure of $\sigma (E)$)
be a curve. Note that $\theta : Z \to \G$ factors through
quasi-finite morphisms $Z \to C$ and $C \to \G$. (2) Choosing
``bigger" $F$ and
$\G$ one can always suppose that 
$\G^* , Z^*$, and $S \setminus F$ are smooth. \vs

\noin {\bf 3. Finiteness of $\theta$. } \vs

Let $S$ and $F$ be as in Lemma 2.3. In the case of $S\simeq \C^2$
there is an isomorphism $H_0(F) \simeq H_3(S \setminus F)$: namely
we assign to each $x_0 \in F$ the 3-cycle presented by a small
sphere with center at $x_0$. 
The same remains true in the normal case.

\lemma 3.1. {\em Let $S$ be a normal affine surface, $\G$ be a
closed curve in $S$, $ \breve{S} =S\setminus \G$, $F$ be a finite
subset of $\G$ such that $\G^* = \G \setminus F$  and
$\tS = S \setminus F$ are smooth. 

{\rm (i)} There exist a natural isomorphism $\phi : H_0(F) \to H_3
(\tS )$.

{\rm (ii)} Let $\varphi : H_3( \tS ) \to H_1(\G^*)$ be the
composition of the Thom isomorphism $H_3(\tS ,  \breve{S}) \to
H_1( \G^* )$ and the natural homomorphism $H_3(\tS ) \to H_3(\tS ,
\breve{S})$, let $x_0 \in F$, $\omega = \phi (x_0)$, and $\G_j, \,
j=1, \ldots , k$ be irreducible analytic branches of $\G$ at
$x_0$. Suppose that $\gamma_j$ is a simple loop in $\G_j$ around
$x_0$ with positive orientation. Then $\varphi (\omega )=
\sum_{j=1}^k \gamma_j$ (where we treat $\gamma_j$ as an element of
$H_1(\G^* )$).}

\proof As $H_3( S)=H_4(S)=0$ [Mil, Th. 7.1] 
the exact homology sequence of pair $(S, \tS )$ shows that
$H_3(\tS ) \simeq H_4(S, \tS )$. Thus in  (i) it
suffices to prove that $H_4(S, \tS ) \simeq H_0(F)$. 
As $S$ may not be smooth we cannot apply Thom's isomorphism. But
by the excision theorem $H_4(S,\tS ) \simeq H_4(U, U\setminus F)$
where $U$ is a small Euclidean neighborhood of $F$ in $S$, i.e.
$U$ is a disjoint union $\bigcup_i U_i$ of neighborhoods $U_i$ of
a point $x_i$ running over $F$. Hence $H_4(S,\tS ) = \oplus_i
H_4(U_i, U_i\setminus x_i)$.  Let $B_{\epsilon}$ be the ball of
radius $\epsilon$ in $\C^n \supset S$ with center at $x_0 \in F$
and $U_0 = S \cap B_{\epsilon}$. For small $\epsilon >0$, $\omega
= S \cap \de B_{\epsilon}$ is a three-dimensional real manifold
and $U_0$ is a cone over $\omega$ [Lo]. Hence $H_4(U_0 , U_0
\setminus x_0) \simeq H_3(U_0 \setminus x_0) \simeq \Z$ and $U_0$
and $\omega$ are generators of $H_4(U_0 , U_0 \setminus x_0)$ and
$H_3(U_0 \setminus x_0)$ respectively. Thus $H_4(S,\tS ) \simeq
H_0 (F)$ which is (i). For small $\epsilon >0$ the sphere $\de
B_{\epsilon}$ and, thus, $\omega$ meet each $\G_j$ transversally
along a simple loop $\gamma_j$ around $x_0$ in $\G_j$. Hence the
description of Thom's isomorphism [Do] yields $\varphi (\omega
)=\sum_{j=1}^k \gamma_j$. \qed

\prop 3.2. {\em Let the assumption and notation of Lemma 2.3 hold,
$H_3(X')=0$, and $\tX =\pi^{-1}( \tS )$ be smooth where $\tS =
S\setminus F$.

{\rm (i)} Then $\theta$ is a finite morphism.

{\rm (ii)} If $x_0$ is a selfintersection point of $\G$ then
$\theta^{-1}(x_0)$ consists of one point $z_0$.


{\rm (iii)} Let $X'$ be also smooth and $H_2(X')=0$. Let $Z^1$ be
an irreducible component of $Z$ and $\G^1 = \theta (Z^1)$. If
$\theta |_{Z^1} : Z^1 \to \G^1$ is not bijective then $\G^1$ is a
polynomial curve (i.e. its normalization is $\C$) and $Z^1$ is
rational.}


\proof Put $ \breve{X} = X' \setminus E$ and $
\breve{S}=S\setminus \G$, i.e. $ \breve{X } \simeq \breve{S}
\times \C$. 
We have the following commutative diagram of the exact homology
sequences of pairs.

\begin{picture}(200,95)

\put(0,65){$\dots\longrightarrow H_{j+1}(\tX, \breve{X})
\longrightarrow H_{j}(\breve{X})\longrightarrow H_{j}(\tX)
\longrightarrow H_{j}(\tX, \breve{X}) \longrightarrow
H_{j-1}(\breve{X})\longrightarrow \dots$}

\put(63,50){$\vector(0,-1){25}$} \put(137,50){$\vector(0,-1){25}$}
\put(199,50){$\vector(0,-1){25}$} \put(122,35){$\simeq$}
\put(270,50){$\vector(0,-1){25}$}
\put(342,50){$\vector(0,-1){25}$} \put(327,35){$\simeq$}

\put(5,5){$\dots\longrightarrow H_{j+1}(\tS , \breve{S})
\longrightarrow H_{j}(\breve{S})\longrightarrow \, H_{j}(\tS ) \,
\longrightarrow \, H_{j}(\tS , \breve{S}) \, \longrightarrow \,
H_{j-1}(\breve{S}) \, \longrightarrow \dots$}

\end{picture}

\noindent which we consider for $j=3$. By Remark 2.4 we can
suppose that $\G^*$, $\tS$, and $E^*$ are smooth. As $\tX
\setminus \breve{X}= E^*$ the Thom isomorphism [Do] implies $H_i(
\tX , \breve{X}) \simeq H_{i-2}(E^*)$ and, similarly, $H_i( \tS ,
\breve{S}) \simeq H_{i-2}(\G^*)$. As $H_3( \breve{S})=0$ [Mil, Th.
7.1] and $H_i(E^*)=H_i(Z^*)$  we have

\begin{picture}(200,95)

\put(10,65){$\dots \longrightarrow \, \, \, 0 \, \, \,
\longrightarrow H_{3}(\tX ) \longrightarrow
H_{1}(Z^*)\longrightarrow H_{2}(\breve{X})\longrightarrow H_2(\tX
)\longrightarrow \dots$} \put (141,
75){$\varphi'$}\put(201,75){$\psi'$}\put(263,75){$\chi'$}
\put(60,50){$\vector(0,-1){25}$} \put(109,50){$\vector(0,-1){25}$}
\put(114,35){$\delta_3$}\put(170,50){$\vector(0,-1){25}$}
\put(175,35){$\theta_1$} \put(235,50){$\vector(0,-1){25}$}
\put(220,35){$\simeq$} \put(240,35){$\breve{\delta}_2$}
\put(295,50){$\vector(0,-1){25}$}\put(300,35){$\delta_2$}
\put(400,35){$(*)$} \put(10,5){$\dots \longrightarrow \, \, \, 0
\, \, \, \longrightarrow \, H_{3}(\tS ) \, \longrightarrow \,
H_{1}(\G^*) \longrightarrow \, \, H_{2}(\breve{S})\longrightarrow
 H_2(\tS ) \longrightarrow
\dots$}\put(142,15){$\varphi$}\put(204,15){$\psi$}
\put(266,15){$\chi$}
\end{picture}

\noindent As $H_4(X')=0$ [Mil, Th. 7.1] the exact homology
sequence of pair $(X',\tX )$ implies that $H_3(\tX )$ is
isomorphic to $H_4(X',\tX )$. By Lemma 2.1 $L= \pi^{-1}(F) = E
\setminus E^*=X' \setminus \tX$ is a curve. If $U$ is a small
Euclidean neighborhood of $L$ in $X'$ then $H_4(U,U\setminus
L)=H_4(X',\tX )$ by the excision theorem. That is, every 3-cycle
$\omega' \in H_3(\tX )$ can be chosen in $U\setminus L$. Thom's
isomorphism maps the image of $\omega'$ in $H(U, U\setminus L)$ to
a 1-cycle $\alpha \in H_1(E^* \cap U)$. As Thom's isomorphisms are
functorial under open embeddings [Do, Ch. 8, 11.5] we can treat
$\alpha$ as an element of $H_1(E^*)$. Thus the image of $\alpha$
under the homomorphism generated by the natural projection $E^*
\to Z^*$ coincides with $\beta_0 = \varphi' (\omega' ) \in
H_1(Z^*)$. As $\alpha$ is contained in a small neighborhood $U
\cap E$ of $L$ in $E$ we see that $\beta_0$ is contained in a
small neighborhood of the finite set $Z \setminus Z^*$. This
implies that ${\rm Im} \, \varphi'$ is contained in the kernel
$G'$ of the homomorphism $H_1(Z^*) \to H_1(Z)$ induced by the
natural embedding $Z^* \hookrightarrow Z$. Let $ \bar{\theta} :
\bar{Z} \to \bar{\G}$ be the extension of $\theta$ to completions
of $Z$ and $\G$. Assume that $F' = (\bar{Z}\setminus Z) \setminus
\bar{\theta}^{-1}( \bar{\G} \setminus \G )$ is not empty which is
equivalent to the fact that $\theta$ is not finite. 
Let $z_0 \in F'$. One can suppose that $x_0 = \theta (z_0) \in F$.
Let $\G_j , \, j=1, \ldots , k$ be irreducible analytic branches
of $\G$ at $x_0$, and let $\gamma_j$ be a simple loop in $\G_j$
around $x_0$ with positive orientation. We treat $\gamma :=
\sum_{j=1}^k \gamma_j$ as an element of $H_1(\G^* )$. 
By Lemma 3.1 $H_3(\tS )$ is generated by elements $\omega$ of form
$\omega= \phi (x_0)$ and $\varphi (\omega )= \gamma$. Take a
simple loop $\beta_1$ around $z_0$ in an irreducible analytic
branch $Z_1$ of $ \bar{Z}$ at $z_0$ such that $ \bar{\theta}
(Z_1)$ is, say, $\G_1$. Take other simple loops around the points
of $ \bar{\theta}^{-1}(x_0)$ whose images under $\theta$ are
contained in $\bigcup_{j=2}^k \G_j$. Then we can construct $\beta
\notin G'$ as an integer linear combination of $\beta_1$ and these
other loops so that $\theta_1 (\beta) =m \gamma$ where $m>0$. As
$\gamma \in {\rm Im} \, \varphi$ we have $\psi (\gamma) =0$. As $
\breve{\delta}_2$ is an isomorphism we have $\psi'(\beta)=0$, i.e.
$\beta \in {\rm Im} \, \varphi'$. This contradiction implies (i).

For (ii) we need another look at ${\rm Im} \, \varphi'$. Let
$Z\setminus Z^* = \{ z_i \}$ and $G_i'$ be the kernel of the
homomorphism $H_1(Z \setminus z_i) \to H_1(Z)$ induced by the
natural embedding $Z\setminus z_i \hookrightarrow Z$. Note that
$G'$ is naturally isomorphic to $\bigoplus_i G_i'$. Moreover,
since $L$ is the disjoint union $\bigcup_i \vartheta^{-1}(z_i)$
one can choose a neighborhood $U$ of $L$ in $X'$ as a disjoint
union of neighborhoods of the curves $\vartheta^{-1}(z_i)$. Hence,
repeating the argument with the excision theorem, we get ${\rm Im}
\, \varphi' = \oplus_i (G_i' \cap {\rm Im} \, \varphi')$. Let
$x_0$ be a selfintersection point of $\G$ (i.e. the number $k$ of
irreducible analytic branches $\G_j$ of $\G$ at $x_0$ is at least
2) and assume for simplicity that $\theta^{-1} (x_0)$ consists of
two points $z_0$ and $z_1$. Note that $z_0,z_1 \in Z \setminus
Z^*$. Let $\beta' \in G_0' \cap {\rm Im} \, \varphi'$ (i.e.
$\beta' =\varphi' (\omega' )$). Then $\delta_3(\omega' )=m\omega$
where $m \in \Z$ and $\omega$ is as before. As $\varphi (m\omega)
= \theta_1 (\beta')$ we see that $\theta_1(\beta') =m \gamma$.
Similarly $\theta_1(\beta'' )$ is proportional to $\gamma$ for
every $\beta'' \in G_1' \cap {\rm Im} \, \varphi'$. On the other
hand since $x_0$ is a selfintersection point, we can find $\gamma'
\in G_0'$ and $\gamma'' \in G_1'$ such that $\theta_1 (\gamma' )$
and $\theta_1 (\gamma'' )$ are not proportional to $\gamma$ but
$\theta_1 (\gamma' +\gamma'')$ is. The commutativity of (*)
implies that $\gamma' + \gamma'' \in {\rm Ker} \, \psi' = {\rm Im}
\, \varphi'$. But the decomposition of ${\rm Im} \, \varphi'$ we
discussed before, implies that $\gamma'$ and $\gamma''$ must be in
${\rm Im} \varphi'$. This contradiction implies (ii).


%
Let $G$ be the kernel of the homomorphism $H_1( \G^* ) \to H_1 (\G
)$, induced by the natural embedding. There is a non-canonical
isomorphism $H_1(\G^* ) = G\oplus H_1(\G^{\n})$ (resp. $H_1( Z^*)=
G' \oplus H_1( Z^{\n})$) such that $\Im \varphi$ (resp. $\Im
\varphi'$) is contained in the first term and $\theta_1 (G')$ is a
subgroup of finite index in $G$. As $\breve{\delta}_2$ is an
isomorphism the ranks of $H_1( \G^{\n})$ and $H_1(Z^{\n})$ are the
same, and the ranks of $G/ \Im \varphi $ and $G' / \Im \varphi' $
coincide as well. For (iii) note that (i) and (ii) implies that
$\theta |_{Z^1} : Z^1 \to \G^1$ is bijective as soon as it is
birational. If it is not birational then the fact, that the ranks
of $H_1( \G^{\n})$ and $H_1(Z^{\n})$ are the same, implies that
the genus of $Z^1$ (and, therefore, $\G^1$) is zero. Assume that
the normalization of $\G^1$ is different from $\C$, i.e. one can
suppose that $ \bar{\G}^1 \setminus \G^1$ consists of more than
one point where $ \bar{\G}^1$ is the closure of $\G^1$ in $
\bar{\G}$. Let $ \bar{Z}^1$ be the closure of $Z^1$ in $ \bar{Z}$.
Then either $ \bar{ \theta }$ maps $\bar{Z}^1 \setminus Z^1$
bijectively onto $\bar{\G}^1 \setminus \G^1$ or it does not. In
the first case a simple loop in $\G^1$ around a point $x \in
\bar{\G}^1 \setminus \G^1$ produces $\gamma^1 \in H_1(\G^*)$ which
is not in $G + {\rm Im} \, \theta_1$ and in the second case an
integer linear combination of simple loops in $Z^1$ around points
of $ \bar{\theta}^{-1}(x)$ produces $\beta^1 \in H_1(Z^*)$ which
is in ${\rm Ker} \, \theta_1 \setminus G'$. In the first case
since $G \supset {\rm Im} \, \varphi ={\rm Ker} \, \psi$, we see
that $\varphi (\gamma^1 )$ is not in ${\rm Im }\, \breve{\delta}_2
\circ \psi' $. On the other hand, as $X'$ is smooth and
$H_2(X')=0$ we get $H_2(\tX )=0$ whence ${\rm Im} \, \psi' = H_2 (
\breve{X})$. This is a contradiction since $ \breve{\delta}_2$ is
an isomorphism. The second case is similar which concludes (iii).
\qed

%
%
%

\rem 3.3. Proposition 3.2 implies that $\pi$ is surjective which
is a generalization of the result of Bonnet [Bo] who proved
this in the case of $X' \simeq \C^3$. 
This implies that when $X'$ is Cohen-Macaulay,  the algebra of
regular functions on $X'$ is faithfully flat over the algebra of
regular functions on $S$ (this follows from [Ma, Chap. 2, (3.J)
and Th. 3] and [Ei, Th. 18.16], see also [Da]). \vs

%

\noin {\bf 4. Injectivity of $\theta$. \\} \lemma 4.1. {\em Let
the notation of Lemma 2.3 hold, $\cM \subset E$ be a germ of a
curve at a general point $x'$ of $E$ such that $\pi |_{\cM} : \cM
\to S$ is injective. Let $\cP \subset X'$ be a germ of a surface
at $x'$ that meets $E$ transversally along $\cM$. Then $\pi
|_{\cP} : \cP \to S$ is injective.}

\proof While speaking of affine modifications here and in Section
5 we use the terminology from [Ka02]. Consider the affine
modification $\sigma : X' \to X \simeq
S \times \C$ from Lemma 2.2. 
Proposition 2.12 from [Ka02] implies that in a neighborhood of a
general point of $E$ this modification $\sigma$ is a composition
of basic modifications. 
Locally, these basic modifications are nothing but usual monoidal
transformations [Ka02, Remark 2.3] (see also Remark 5.2 below). In
the case of three-dimensional manifolds and one-dimensional
centers the coordinate form of these monoidal transformations is
$(\xi_1, \xi_2 , \xi_3) \to (\xi_1, \xi_1 \xi_2 , \xi_3) : =(
\zeta_1,
\zeta_2, \zeta_3)$. 
If $\sigma$ is locally such a monoidal transformation, then
$\sigma$ maps $\cM$ injectively into the line $\zeta_1 = \zeta_2
=0$ since $\pi = {\rm pr} \circ \sigma$ where ${\rm pr} : X \to S$
is the natural projection. Hence $\cM$ must be given by equations
of form $\xi_1= \xi_2 -g_0 (\xi_3)=0$. As $\cP$ meets $E$
transversally along $\cM$ we see that $\cP$ is given by an
equation of form $\xi_2 -g(\xi_1 , \xi_3)$ where
$g(0,\xi_3)=g_0(\xi_3)$. Hence the image of $\cP$ is the germ of a
surface given by $\zeta_2 =\zeta_1 g(\zeta_1, \zeta_3)$. In
particular, this germ is isomorphic to $\cP$ (via the
transformation) and transversal to the plane $\zeta_1=0$.
Induction on the number of monoidal transformations in the local
decomposition of $\sigma$ implies that $\sigma |_{\cP}: \cP \to
\cQ := \sigma (\cP )$ is bijective and $\cQ$ meets $D$
transversally along a germ $\cC$ of $C=\sigma (E)$ at $x= \sigma
(x')$. As $x$ is a general point of $C$, projection ${\rm pr}$
maps $\cC$ injectively into $\G$. As $\cQ$ meets $D$ transversally
along $\cC$, morphism ${\rm pr}|_{\cQ} : \cQ \to S$ is not
ramified along $\cC$ whence this map is injective. Thus $\pi
|_{\cP} : \cP \to S$ is injective. \qed

\prop 4.2. {\em Let the assumption of Proposition 3.2 (iii) hold,
$S$ be smooth, and $\pi : X' \to S$ be the quotient morphism of a
free $\C_+$-action $\Phi$ on $X'$. Then $\theta$ is bijective.}

\proof Assume the contrary, i.e. there is a component $Z^1$ of the
$Z$ such that $\theta |_{Z^1}$ is not bijective. By Proposition
3.2 (i)-(ii) $\theta |_{Z^1}$ is not birational. Furthermore,
Proposition 3.2 (iii) implies that $Z$ contains a ramification
point $z_0$ of $\theta$ whose order is $m\geq 2$. Take a smooth
point $x_0'$ of $\vartheta^{-1}(z_0)$. The irreducible component
$L_0$ of $\vartheta^{-1}(z_0)$, that contains $x_0'$, is an
irreducible component of a fiber of $\pi$ and, therefore, an orbit
of $\Phi$ since $\Phi$ is free. Consider the germ $ {\cal P}$ of a
surface at $x_0'$ transversal to $L_0$. Then all orbits of $\Phi$
close to $L_0$ must be transversal to $\cP$. Let us show that this
cannot be true.

Consider $\varrho = \pi |_{{\cal P}} : {\cal P} \to \pi ({\cal P}
)$. One can see that $\varrho$ is finite (as $\cP$ does not
contain fibers of $\pi$). 
As $\pi |_E = \theta \circ \vartheta$ and $\theta$ has a
ramification point of order $m$ we see that $\varrho$ is at least
$m$-sheeted. As $S$ is smooth there must be a ramification curve
$R \subset \cP$ of morphism $\varrho$. Consider a general point
$x_1'$ of the curve $\cM=\cP \cap E$, and the germ $\cP_1$ of
$\cP$ at $x_1'$. Note that $\cP_1$ meets $E$ transversally along
the irreducible germ $\cM_1$ of $\cM$ at $x_1'$, and $\pi
|_{\cM_1}$ is injective. By Lemma 4.1, $\varrho$ is locally
injective in a neighborhood of $x_1'$, i.e. $R$ has no common
components with ${\cal M}$. Therefore, a general point $x_2' \in
R$ is contained in $X'\setminus E \simeq \breve{S} \times \C$. As
$x_2'$ is a ramification point of $\varrho$ it cannot be a point
where $\cP$ meets the fiber of the natural projection $ X'
\setminus E \to \breve{S}$ (which is a fiber of $\pi$ and,
therefore, an orbit of $\Phi$) transversally which yields the
desired contradiction. \qed

\noin {\bf 5. Main theorem. \\} We want to decompose $\sigma : X'
\to X$ from Lemma 2.2 into a composition of simpler modifications
$X' \to X_1$ and $X_1 \to X$. Let $A$, $A_1$, and $A'$ be the
algebras of regular functions on $X$, $X_1$, and $X'$
respectively. In order to find $X_1$ it suffices to find an affine
domain $A_1$ which is a proper subring of $A'$ and which contains
but not equal to $A$.

\lemma 5.1. {\em Let $\sigma : X' \to X$ be an affine
modification, $C$ be its center, $D^1$ be an irreducible component
of its divisor $D$, $C^1 =C \cap D^1$. Suppose that $A' =A[I/f]$
where $I$ is an ideal in $A$ that vanish on $C$, $f=g^{n} h$, $D
=f^{-1}(0)$, and $D^1=g^{*}(0)$ while $h$ is not identically zero
on $D^1$. Let $E^1 = \sigma^{-1} (C^1)$, and the defining ideal
$\cI_{A'}(E^1)$ of $E^1$ be the principal ideal in $A'$ generated
by $g$.

{\rm (1)} Then $A[(\cI_A(C^1))^{n}/g^{n}] \subset A'$ for the
defining ideal $\cI_A(C^1)$ of $C^1$ in $A$. In particular, for
every ideal $J \subset \cI_A (C^1)$ the algebra $A_1 =A[J/g]$ is
contained in $A'$.

{\rm (2)} If $J$ in (1) coincides with $\cI_A (C^1)$ then $A'
=A_1[I_1/f_1]$ where $I_1$ is an ideal in $A_1$ and $f_1 =
g^{n-1}h$.}

\proof Note that for every $a \in (\cI_A(C^1)) ^{n}$ we have $a
\in (\cI_{A'}(E^1)) ^{n}$ (as an element of $A'$) whence $a/g^{n}
\in A'$. For the second statement of (1) note that $A_1
=A[(g^{n-1}J)/g^{n}] \subset A[(\cI_A(C^1))^{n}/g^{n}]$.

Let $g, b_1, \ldots , b_k$ be generators of $I \subset
\cI_A(C^1)$. Then in (2) $b_1/g, \ldots , b_k/g \in A_1$. The
ideal $I_1$ generated by these elements is what we need in (2).
\qed

\rem 5.2. Let $x \in C^1$ be a smooth point of both $C^1$ and $X$,
and the zero multiplicity of $g$ at $x$ be 1. In a small Euclidean
neighborhood $U$ of $X$ the defining ideal of $C^1 \cap U$ is a
complete intersection, i.e. it is generated by regular functions
$g, b_1, \ldots , b_k$ which are elements of a local coordinate
system on $U$. Hence in $U$ the modification $A \hookrightarrow
A_1=A[\cI_A (C^1)/g]$ can be viewed as a monoidal transformation.
This fact in combination with Lemma 5.1 (2) and induction on $n$
yields another way to show that modification $\sigma$ in the proof
of Lemma 4.1 is a composition of monoidal transformations in a
neighborhood of a general point of its center. \vs

\prop 5.3. {\em Let the assumption of Proposition 3.2 (iii) hold,
and $\G^1$ (resp. $Z^1$) be an irreducible components of $\G$
(resp. $Z$) such that $\theta ( Z^1) \subset \G^1$. Suppose that
$\theta |_{Z^1} : Z^1 \to \G^1$ is bijective. 
Let $\G^{\bullet}$ be the closure of $\G \setminus \G^1$ in $S$.
Then $\pi^{-1}(S \setminus \G^{\bullet})$ is isomorphic to $(S
\setminus \G^{\bullet}) \times \C$ over $S \setminus
\G^{\bullet}$.}

\proof Let $\sigma : X' \to X \simeq_S S \times \C$ be the affine
modification from Lemma 2.2. In particular, its divisor $D=\G
\times \C$ contains an irreducible component $D^1 = \G^1 \times
\C$. Our aim is to construct another affine modification
$\sigma^{\bullet} : X' \to X$ such that its divisor is
$D^{\bullet} = \G^{\bullet} \times \C$ which would yield the
statement of this Proposition. As $S$ is factorial by Lemma 2.1,
$\G^1$ is the zero locus of a regular function $g$ on $S$. If we
treat $g$ as a regular function on $X$ we get $A' = A[I/f]$ where
$I$ and $f=g^{n} h$ are as in Lemma 5.1. In virtue of [Ka02, Lemma
2.2] one can suppose that the center $C$ of $\sigma$ is a curve.
By Remark 2.4 $\theta : Z \to \G$ factors through ${\rm pr}|_C : C
\to \G$ and $Z \to C$. Let $C^1$ be the image of $Z^1$ under the
last morphism. Then $\pr |_{C^1}: C^1 \to \G^1$ is bijective.
Hence $C^1$ is given in $D^1= \G^1 \times \C$ by $t-{\bar{b}}_1=0$
where $t$ is a coordinate on the second factor of $X\simeq S\times
\C_t$ and $ \bar{b}_1$ is a continuous rational function on
$\G^1$. Let
us show that $t-\bar{b}_1$ is regular on $D^1$. 
Let $x_0$ be the preimage of $o \in \G^1$ in $C^1$ and $x_0' \in
\pi^{-1}(o)$ be a smooth point of $\pi^{-1}(o)$. Take the germ
$\cP$ of an analytic surface at $x_0'$ such that $\cP \cap L
=x_0'$. One can check that $\cQ = \sigma (\cP )$ is closed in a
neighborhood of $x_0$ and meets $D_1$ along the germ $\cC^1$ of
$C^1$ at $x_0$. Hence $\cQ$ is given in $X\simeq S \times \C_t$ by
zeros of a function $t^n+a_{n-1}t^{n-1} + \ldots + a_0 $ where
each $a_i$ is a germ of a continuous function at $o \in S$ which
is analytic outside $\G^1$. As $S$ is normal $a_i$ is analytic at
$o$ as well. Hence by [Ka02, Lemma 3.3] $t- \bar{b}_1$ is regular.
Let $b_1$ be an extension of $ \bar{b}_1$ to $S$. Then
$A_1=A[(t-b_1)/g]$ is the algebra of regular functions on a
hypersurface $X_1$ in $X \times \C_u$ given by equation $ug
=t-b_1$, and 
$X_1 \simeq S \times \C_u \simeq X$. By Lemma 2.1 the defining
ideal of $E^1=\sigma^{-1}(C^1)$ in $A'$ is the principal ideal
generated by $g$. As Lemma 2.3 implies that $Z^1=\theta^{-1}(\G^1
)$, we have $C^1=D^1\cap C$. Thus by Lemma 5.1 $A_1 \subset A'$.
Hence we can replace modification $\sigma : X' \to X$ by
modification $X' \to X_1\simeq_S S\times \C$ (induced by the
embedding $A_1 \hookrightarrow A'$) and all assumption of
Proposition 3.2 hold but the advantage is that now $f=g^{n-1} h$
by Lemma 5.1 (since the defining ideal of $C^1$ in $A$ is
generated by $g$ and $t - b_1$). Repeating this procedure we can
make $n=0$, i.e. we get the modification $\sigma^{\bullet} : X'
\to X$ which concludes the proof. \qed

\thm 5.4. {\em Let the assumption of Lemma 2.3 hold, $X'$ be
smooth, and $H_2(X')=H_3(X')=0$.

{\rm (i)} Then $\G$ can be chosen so that each irreducible
component $Z^1$ of $Z$ (resp. $\G^1 =\theta (Z^1)$ of $\G$) is
rational (resp. a polynomial curve), and $\theta |_{Z^1} :Z^1 \to
\G^1$ is not bijective.

{\rm (ii)} If $S$ is smooth and the $\C_+$-action $\Phi$ on $X'$
has no fixed points then $X'$ is isomorphic to $S\times \C$ over
$S$ and the action is generated by a translation on the second
factor.}


\proof Applying Proposition 5.2 consequently to each component
$Z^1$ of $Z$, such that $\theta |_{Z^1} : Z^1 \to \G^1 =\theta
(Z^1)$ is bijective, we obtain (i) in virtue of Proposition 3.2
(iii). Proposition 4.2 and (i) imply that the quotient morphism
$\pi : X' \to S$ makes $X'$ isomorphic
to $S \times \C$ over $S$. \qed 

{\bf Proof of Theorem 1.} By [Miy80] $\C^3//\C_+ \simeq \C^2$ for
every nontrivial $\C_+$-action on $\C^3$ whence Theorem 5.4 (ii)
implies the desired conclusion. \qed

\rem 5.5. (1) There is another way to extract Theorem 1 from
Proposition 4.2. For a bijective $\theta$ and $X'\simeq \C^3$ an
easy Euler characteristics argument shows that all fibers of
$\vartheta$ (not only the general ones) are irreducible. This, in
turn, implies that the $\C_+$-action is separated and, therefore,
it is a translation [DeFi, Kr].

(2) Consider an $n$-dimensional smooth contractible affine
algebraic variety $X'$ and a free algebraic action of a unipotent
group $U$ on it (i.e. each orbit of the action is isomorphic to
$U$). Suppose that $U$ is of dimension $n-2$ (i.e. $U$ is
isomorphic to $\C^{n-2}$ as an affine algebraic variety). By
Fujita's result (e.g., see [Ka94, Prop. 3.2]) $X'$ is factorial
and by Zariski's theorem [Za] $X'//U$ is still an affine algebraic
surface whence one can repeat all Lemmas and Propositions before
for this situation with obvious adjustments (say, in Lemma 2.3
$E^*$ will be isomorphic to $Z^* \times \C^{n-2}$). As Miyanishi's
theorem [Miy80] says that in the case of a free action $\C^n//U
\simeq \C^2$, we have the following analogue of Theorem 1: for
every free algebraic action of $U$ on $\C^n$ there is an
isomorphism $\C^n \simeq \C^2 \times U$ such that the action is
induced by the natural action on the second factor.

\begin{center} REFERENCES \end{center} \vs
\noin [Ba] H. Bass, {\em A nontriangular action of $G\sb{a}$ on
$A\sp{3}$}, J. Pure Appl. Algebra {\bf 33} (1984), no. 1, 1--5. \\
\noin [Bo] P. Bonnet, {\em Surjectivity of quotient maps for
algebraic $(\C ,+)$-actions and
polynomial maps with contractible fibers}, Trans. Groups (to appear). \\
\noin [Da] D. Daigle, {\em A note one generic planes and locally
nilpotent derivations of $k[X,Y,Z]$}, (to appear). \\
\noin [DaFr] D. Daigle, G. Freudenburg, {\em Locally nilpotent
derivations over a UFD and an application to rank two locally
nilpotent derivations of $k[X_1,\ldots ,X_n]$},
J. Algebra {\bf 204} (1998), no. 2, 353--371. \\
\noin [DeFi] J. Deveney, D. Finston, {\em Free $G_a$ actions on
$\C^3$}, Proc. Amer.
Math. Soc. {\bf 128} (2000), no. 1, 31--38. \\
\noin [Do]{Do} A. Dold, {\em Lectures on Algebraic Topology},
Springer, Berlin e.a., 1974. \\
\noin [Ei] D. Eisenbud, {\em Commutative Algebra with a View
Towards Algebraic Geometry}, Graduate Texts in Math., Springer,
N.Y. e.a., 1994. \\
\noin [Ka94] S. Kaliman, {\em Exotic analytic structures and
Eisenman intrinsic measures}, Israel Math. J. {\bf 88} (1994),
411--423.\\
\noin [Ka02] S. Kaliman, {\em Polynomials with general
$\C^2$-fibers are variables}, Pacific J. Math., {\bf 203} (2002), no. 1, p. 161--189.\\
\noin [Kr] H. Kraft, {\em Challenging problems on affine
$n$-space}, S\'eminaire
Bourbaki, Vol. 1994/95; Ast\'erisque {\bf 237} (1996), 295--317. \\
\noin [Lo] E.J.N. Looijenga, {\em Isolated singular
points on complete intersections}, London Mathematical Society
Lecture Note
Series, {\bf 77}, Cambridge University Press, Cambridge, 1984, 200 p. \\
\noin [Ma] H. Matsumura, {\em  Commutative Algebra}, The Benjamin/Cummings Publisthing
Company, Reading, Massachusetts, 1980.\\
\noin [Mil] J. Milnor, {\em Morse Theory}, Princeton Univ. Press,
Princeton, NJ, 1963. \\
\noin [Miy80] M. Miyanishi, {\em Regular subrings of a polynomial
ring}, Osaka J. Math. {\bf 17} (1980), no. 2, 329--338. \\
\noin [Re] R. Rentschler, {\em Op\'erations du groupe additif sur
le plane affine}, C.R. Acad. Sci. Paris, {\bf 267} (1968),
384--387.\\
\noin [Sn] D. Snow, {\em Unipotent actions on affine space}, in
``Topological methods in algebraic transformation groups",
Birkh\"auser, Boston, MA, 1989, 165--176. \\
\noin [Za] O. Zariski, {\em Interpretations
alg\'ebrico-g\'eom\'etriques du quatorzi\`eme probl\`eme de
Hilbert}, in ``O. Zariski: Collected Papers", {\bf 2}, The MIT
Press, Cambridge (Massachusetts), 1973, 261--275.

\end{document}